\documentclass[12pt]{article}
\usepackage{graphicx}
\usepackage{harvard}
\usepackage{array}
\usepackage{fullpage}
\usepackage[small, hang, bf]{caption}
\newcommand{\Keywords}[1]{\par\noindent{\small{\em Keywords\/}: #1}}

\begin{document}

\title{Estimating Omissions from Searches}

\author{
Anthony J. Webster\thanks{email: anthony.webster@ccfe.ac.uk} 
{ }and 
Richard Kemp\\ 
\normalsize{United Kingdom Atomic Energy Authority,}\\
\normalsize{Culham Science Centre, Abingdon, Oxon, OX14 3DB.}
}

\date{Published in:\\ The American Statistician, \\Volume 67, Issue 2,
  pages 82-89, 2013.}

\maketitle

\begin{abstract}

The mark-recapture method was devised by Petersen in 1896 to estimate
the number of fish migrating into the Limfjord, and independently by
Lincoln in 1930 to estimate waterfowl abundance. The technique can be
applied 
to any search for a finite number of items by two or more people or
agents, allowing the number of searched-for items to be estimated.  
This ubiquitous problem appears in fields from ecology and
epidemiology, through to mathematics, social sciences, and
computing. Here we exactly calculate the moments of the hypergeometric
distribution associated with this long-standing problem, 
confirming that
widely used estimates conjectured in 1951 are often too small. Our
Bayesian approach highlights how different search strategies will
modify the estimates. 
The estimates are applied to several examples.  
For some published applications substantial errors are found to result 
from using the Chapman or Lincoln-Petersen estimates. 

\vspace{0.5cm}


\vspace{0.5cm}

\end{abstract}

\section{Introduction}

If a finite set is searched by two or more people it is possible to
estimate how many of the searched-for items have been missed. The
simple Lincoln-Petersen estimate was independently developed by
\citeasnoun{Petersen} to estimate fish numbers migrating between the
German sea and the Limfjord, and by \citeasnoun{Lincoln} to estimate
waterfowl abundance. The technique has rapidly grown in popularity
since a more rigorous treatment by \citeasnoun{1951}, especially in
the context of ecological census techniques \cite{Seber,Sutherland}
and epidemiology \cite{Hook}. Our interest arose from the technique's
application to assess the accuracy of a literature search. In 1938
such a literature search led to the re-discovery of Alexander
Fleming's papers on penicillin \cite{Masters,Lax}, and penicillin's
subsequent development. Today literature searches are a valued method
for identifying and appraising evidence, particularly in
evidence-based healthcare \cite**{EBM}. Reviews often search thousands
of papers, and standardised guidelines have developed for reporting
search terms and the databases used \cite**{Prism,Cochrane}. 
Common
practice involves an electronic search to retrieve hundreds or even
thousands of potentially relevant articles, that are subsequently
searched by the authors for pertinent material. Inevitably, even if
multiple authors search the database, human error may
cause some papers to be erroneously missed at this stage, leading to a
less comprehensive review \cite**{Edwards}. The Lincoln-Petersen
estimator has previously been used to assess the completeness of
medical databases \cite**{Spoor,Bennett,Poorolajal}, and to provide
``stopping rules'' to help determine when searches are complete
\cite**{Kastner,Booth}; 
surprisingly, standard practice does not include an estimate for the
number of papers 
unintentionally omitted by a search.

Here we derive some simple but rigorous results for estimating the
number of items missed from a search, including exact expressions for
the average, standard deviation, and skewness. They correct a widely
used conjecture from \citename{1951}'s 1951 paper and a subsequent
widely used approximation for the variance. Despite their extensive
use \cite{Seber,Hook,Sutherland}, we 
confirm the suggestion \cite{Pelayo} that previous 
conjectured and approximated estimates can  be inaccurate for many
cases of interest, including assessing the accuracy of literature
searches.   

\begin{figure}
\begin{center}
\includegraphics[width=8 cm]{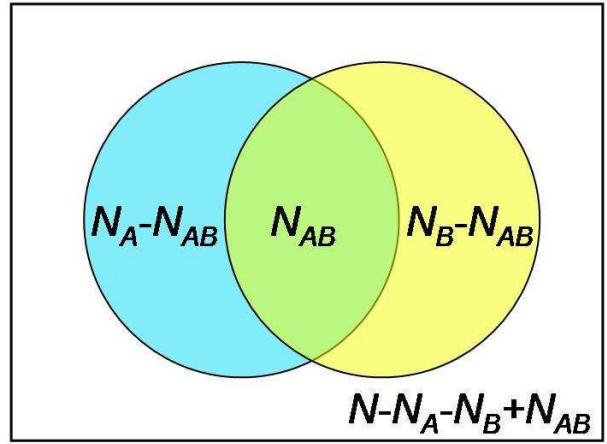}
\end{center}
\caption{ \label{fig:venn1}
The total number of papers found ($N_f$) equals the number found by
$A$ ($N_A$), plus the number found by $B$ ($N_B$), minus the number
found by both $A$ and $B$ ($N_{AB}$), that have been counted twice.  
}
\end{figure}

The problem is as follows. 
Authors A and B each separately search a
given set of references for relevant articles. (It is assumed that
after agreement by both authors, papers that are included are
definitely relevant.) The result is that $N_A$ and $N_B$ articles are
found by authors A and B respectively with $N_{AB}$ of those found by
both authors. If we assume all papers are equally likely to be found,
then a simple estimate can be made as follows. Taking $N$ as the total
number of papers searched for, and taking probabilities $p_A$, $p_B$,
and $p_{AB}$ for 
A, B, and both (A and B) 
finding $N_A$, $N_B$, and $N_{AB}$ papers
respectively, then we can estimate $p_A$, $p_B$, and $p_{AB}$, from 
\begin{equation}\label{probassumpt}
\begin{array}{l}
p_A \approx \frac{N_A}{N} \\
p_B \approx \frac{N_B}{N} \\
p_{AB} \approx \frac{N_{AB}}{N} 
\end{array}
\end{equation}
Because the probability $p_{AB}$ of a paper being found by both
authors is $p_{AB}=p_A \times p_B$, we can combine and solve
(\ref{probassumpt}) for $N$, giving an estimate for $N$ as 
\begin{equation}\label{N1} 
N \approx \frac{N_A N_B}{N_{AB}} 
\end{equation}
The number of papers missed, $X$, is then estimated to be $X=N-N_f$,
where $N_f = N_A+N_B-N_{AB}$ is the total number of different papers
found by both authors (figure \ref{fig:venn1}), finding after a little
algebra, 
\begin{equation}\label{X1}
X \approx \frac{(N_A-N_{AB})(N_B-N_{AB})}{N_{AB}} 
\end{equation}

Equations (\ref{N1}) and (\ref{X1}) are often reasonable estimates if the
numbers involved are large. 
However these estimates are clearly
misleading if $N_{AB}=N_A$, $N_B$, or is zero: for the former cases
because there can be papers that both authors have missed (although
the estimate suggests not); and for the latter case because an
infinite estimate is inconsistent with searching a finite set. 
More importantly,
there is no indication for the accuracy of the estimate, so used in
isolation it is impossible to know whether it is  reasonable or
not. Improved estimates are given later by 
(\ref{av1}), (\ref{sig1}), (\ref{av2}), and (\ref{sig2});
the need for them and their derivation is explained in the following
sections.  
The key assumption underlying all of these estimates is that all items
are equally likely to be found. 
As is discussed at the end of Section \ref{results},  
when this assumption is true or a reasonable approximation, then the estimates
 can be used.

The paper proceeds as follows. 
Section  \ref{sbayes} uses a
Bayesian approach to allow a rigorous mathematical derivation of the
probability density function for the number of items missed.  
Section \ref{results} considers the calculation of its  moments. 
``Exact estimates", refer to exactly calculated moments of the distribution. 
``Approximate estimates", refer to approximations for the moments,   
usually found by expanding about the distribution's maximum. 
Consequently approximated averages are often close to the ``most
probable" estimate, where the distribution is a maximum. 
Section \ref{sbayes2} 
comments on the effects of different assumptions on the final answer,
and finds explicit prior assumptions for which 
Chapman's estimate is exactly the most probable estimate.  
The main result of this paper is to show that the moments can be
calculated exactly, subsequently finding that Chapman's
extensively used estimate can sometimes be misleading. 
A recently published
example discussed in Section \ref{results} emphasises this.

Throughout the paper we refer to two search procedures. 
In the example above, both authors searched for all the papers ($N$) 
and compared the number found by both ($N_{AB}$) to estimate $N\approx
N_AN_B/N_{AB}$.  
An alternative approach is for A and B to search 
for a {\sl predetermined}
number of items $N_A$ and $N_B$ respectively, stopping when that
number is found, 
and again using the number $N_{AB}$ found by both to estimate
$N\approx N_AN_B/N_{AB}$.  
Whereas the former approach is more sensible for a literature search,
the latter approach allows a comparatively small sample of animals to
provide an estimate for their abundance.  
Mathematically the difference can be important. 
If a fixed number of items $N_A$ are searched for, then other than the 
requirement that $N_A \leq N$, $N_A$ is {\sl
  independent} of $N$. In contrast, if all items are searched for then
the probability of A finding $N_A$ items is {\sl dependent} on
$N$. Equivalent remarks apply to B.  
Section \ref{sbayes} uses Bayes theorem to rigorously formulate the
problem for both search procedures.  
Section \ref{results} notes that 
provided that a large number of items are found, then the  
moments of both problems are closely related, and the moments of one
can be used to closely approximate the moments of the other.  
The consequences of different search procedures are discussed
further in Section \ref{sbayes2}.  
Section \ref{conc} summarises the paper's conclusions.

\section{Bayesian formulation}\label{sbayes}

The shortcomings with (\ref{N1}) and (\ref{X1}) arise from the
estimates of $p_A \simeq  N_A/N$, $p_B \simeq  N_B/N$, and $p_{AB}
\simeq N_{AB}/N$. They improve with increasing values of $N_A$, $N_B$,
and $N_{AB}$, but are nonetheless estimates. Specifically, if we know
the probability $p_A$ of author A finding any given paper (we
continue to assume all papers are equally difficult to find), and if
we also knew the total number of papers $N$ that the author is
searching for, then the probability of author A finding $N_A$ papers
is given by the binomial distribution,  
\begin{equation}\label{pna1}
P(N_A | N,p_A) = 
\left( 
\begin{array}{c}
N
\\
N_A 
\end{array}
\right) 
p_A^{N_A} \left( 1 - p_A \right)^{N-N_A} 
\end{equation}
The expected number of papers to be found is then $\langle N_A \rangle
\equiv \sum_{N_A=0}^N N_A P(N_A|N,p_A)=p_A N$ (\emph{e.g.}
\citeasnoun{Stirzaker}). Therefore provided $N_A \simeq \langle N_A
\rangle$, as on {\sl average} it will be, then the estimates
(\ref{probassumpt}) will be reasonable. However, for small
numbers in particular it can give misleading results. 

Bayes' theorem 
was first used for mark and recapture
estimates  by
\citeasnoun{Gaskell}, and allows a rigorous derivation that avoids
these shortcomings. 
In its modern form Bayes' theorem states that $P(X|Y) P(Y)=
P(Y|X)P(X)$ \cite{Sivia}, and allows us to write,  
\begin{equation}
P(N|N_A,N_B, N_{AB}) = \frac{ 
P(N_A,N_B, N_{AB}|N) P(N) }{P(N_A,N_B,N_{AB})} 
\end{equation}
Repeatedly using $P(X,Y)=P(X|Y)P(Y)$ \cite{Sivia}, 
and conditional independence of $N_A$ ($N_A\leq N$), $N_B$ ($N_B \leq
N$), given  $N$,  
this expands to give, 
\begin{equation}\label{pf1}
P(N|N_A,N_B, N_{AB}) 
= 
\frac{P(N_{AB} | N_A,N_B,N) P(N_A|N)P(N_B|N) P(N)}{P(N_A,N_B,N_{AB})} 
\end{equation}
Equation (\ref{pf1}) gives the probability of there being $N$ papers
to find, 
given that author A has found $N_A$ papers, author B has found $N_B$
papers, and $N_{AB}$ of the papers were found by both authors. $P(N)$
is the (prior) probability of there being $N$ papers to be found given
no information about the numbers of papers A and B will find,
$P(N_A|N)$ is the probability of finding $N_A$ papers given that there
are $N$ papers to be found, and equivalently for $P(N_B|N)$. $P(N_{AB}
|N_A,N_B,N)$ is the probability of $N_{AB}$ papers being found by both
authors, given that there are $N$ papers to find, and that authors A
and B each find $N_A$ and $N_B$ papers respectively.

\subsection{Searches for every item}\label{S2.1}

Firstly consider $P(N_A|N)$, and assume that all $N$ items are
searched for. Given  no prior knowledge of how effective 
author A may be at finding papers, we take $P(N_A|N)$ to 
be functionally independent of 
$N_A$. Correct normalisation requires that
$\sum_{N_A=0}^NP(N_A|N)=1$, giving $P(N_A|N)=1/(N+1)$, and similarly
for $P(N_B|N)$. 
Equivalently, 
assume $p_A$ and $N$ are independent, 
and take $P(N_A|N,p_A)$ as given by (\ref{pna1}). 
Then use marginalisation \cite{Sivia} to write 
$P(N_A|N)=\int_0^1 P(N_A|N,p_A)P(p_A) dp_A$,   
assume a uniform prior for $P(p_A)$, and integrate to find the  
same answer. 
This latter approach suggests how 
the method can be generalised if we relax the assumption that all items are
equally likely to be found, through modified forms for $P(N_A|N,p_A)$
and $P(p_A)$.    
$P(N_{AB}|N_A,N_B,N)$ is 
the probability of there being $N_{AB}$ items found by both A and B,
given only the  information that A found $N_A$ items, B found $N_B$
items, and that there are $N$  items to find.    
This can be
calculated by using a metaphor of selecting balls from an urn filled
with $N$ white balls. The first author picks $N_A$ balls at random,
paints them yellow, and returns them. The second author picks $N_B$
balls, and  $N_{AB}$ is the number of yellow balls the second author
has picked. This is a well-known problem (\emph{e.g.}
\citeasnoun[p. 174]{Stirzaker}), whose solution is the hypergeometric
distribution,  
\begin{equation}\label{E20}
P(N_{AB}|N_A,N_B,N) 
= 
\frac{N_A!N_B!(N-N_A)!(N-N_B)!}{N_{AB}!
(N_A-N_{AB})!(N_B-N_{AB})!N!(N-N_f)!} 
\end{equation} 
with $N_{AB} \leq N_A \leq N$ and 
$N_{AB} \leq N_B \leq N$.

Combining the above (\ref{pf1}) and (\ref{E20}) with
$P(N_A|N)=P(N_B|N)= 1/(N+1)$ we get,  
\begin{equation}\label{fulleq}
P(N|N_A,N_B,N_{AB}) 
= 
\frac{ (N-N_A)!(N-N_B)! }{
N!(N-N_f)!} \frac{P(N)}{(N+1)^2} C
\end{equation}
where $C$ is functionally dependent on $N_A$, $N_B$, and $N_{AB}$, but not 
$N$, and is most easily found by ensuring that $P(N|N_A,N_B,N_{AB})$ is
normalised to $1$ after summing over $N$ from the total number of
different papers found $N_f=N_A+N_B-N_{AB}$, to $\infty$. This Bayes'
theory approach was used by \citeasnoun{Zucchini} to derive a similar
result, but without the factors of $P(N_A|N)$ and $P(N_B|N)$ that lead
to some differences discussed later. Note that because the sum is over
$N$ not $N_{AB}$, the moments are different to those usually
associated with the hypergeometric distribution that involve sums over
$N_{AB}$.

\subsection{Searching for a predetermined number of items}

If authors A and B search for a fixed number of say $10$ items each,
so that $N_A$ and $N_B$ are now specified in advance, 
then the previous derivation is modified slightly. 
As before, $N_{AB} \leq N_A \leq N$ and 
$N_{AB} \leq N_B \leq N$, 
but $N$ and $N_{AB}$ can otherwise be 
assumed  independent of $N_A$ and $N_B$. 
If $I$ is some prior information, such as the number of items $N_A$ to
be searched  for by A and the number of items $N_B$ to be searched for
by B, then Bayes' theorem  gives  
 \cite{Sivia} $P(X|Y,I)=P(Y|X,I)P(X|I)/P(Y|I)$. 
Substituting $N$ for $X$, $N_{AB}$ for $Y$, and $N_A,N_B$ for $I$,
Bayes' theorem gives, 
\begin{equation}\label{modBayes}
P(N|N_A,N_B,N_{AB}) 
= \frac{P(N_{AB}|N_A,N_B,N) P(N|N_A,N_B)}{P(N_{AB}|N_A,N_B)} 
\end{equation}
If we make the prior assumption that all values of $N$ (greater than
or equal to 
the largest of $N_A$ and $N_B$), are equally likely, then
$P(N|N_A,N_B)$ will not depend on $N$.  
This is an ``improper'', {\sl i.e.} un-normalisable, prior. 
Strictly $P(N|N_A,N_B)$ should be zero for $N$ bigger than the largest
conceivable number of items in the set being searched. 
With this assumption the  
factor of $P(N|N_A,N_B)$ is replaced with a constant term, leaving,  
\begin{equation}\label{neweq}
P(N|N_A,N_B,N_{AB}) 
= 
\frac{ (N-N_A)!(N-N_B)! }{
N!(N-N_f)!} K
\end{equation}
where, as for $C$ in (\ref{fulleq}), $K$ is functionally dependent on
$N_A$, $N_B$, $N_{AB}$,   and is most easily found by ensuring that
(\ref{neweq}) is correctly 
normalised.  
This is the equation whose 
{\sl approximated}  
moments have been extensively used
\cite{Seber,Sutherland,Hook} and studied
\cite{1951,Zucchini,Seber1970a,Wittes,Pelayo}, and that we will exactly
calculate shortly. 


\section{Results}\label{results}

Given a suitable choice for $P(N)$ or $P(N|N_A,N_B)$ respectively,
(\ref{fulleq}) and (\ref{modBayes}) 
provide the full 
solution to the problem, allowing numerical values for the average and
standard deviation to be calculated by summing from $N=N_f$ to
$N=\infty$ for different moments of $N$. 
The following section takes the prior 
$P(N|N_A,N_B)$ as being constant, then calculates the moments of
(\ref{neweq})  
exactly. 
It also gives an (often excellent) approximation for the moments of
(\ref{fulleq}) when the prior $P(N)$ is constant, and suggests a prior
for which the calculated moments are exact.   
Throughout we will use the statistical physics notation of angled
brackets, with e.g. $\langle f(N) \rangle$, to denote the
expected value of some function $f(N)$, obtained by averaging over the
probability density function for $N$.  
Firstly we will calculate moments of the extensively studied
(\ref{neweq}), and compare these exactly calculated moments with
existing approximations. 
Then we will consider the moments of (\ref{fulleq}), and use these in
some applications. 





\subsection{The moments of (\ref{neweq})}

To calculate the moments we first rewrite (\ref{neweq}) in terms of
$X=N-N_f$, $X_A=N_A-N_{AB}$, and $X_B=N_B-N_{AB}$, so that $N_f =
N_{AB} + X_A + X_B$, and,  
\begin{equation}\label{Papp}
P(X|X_A,X_B,N_{AB}) = 
\frac{(X+X_A)!(X+X_B)!}{X!(X+N_{f})!} 
K
\end{equation}
This gives a probability distribution for the number of papers $X$
that have not been found, with $X$ between $0$ and $\infty$. 
The moments of (\ref{Papp}) are calculated next using a generating
function approach. 
Appendix A 
contains an alternative (our original) calculation
for the moments that is less systematic, but uses
simpler mathematical 
concepts and avoids the use of generating functions. 
All appendices are available as online supplementary material. 
The 
moments of (\ref{Papp}) can be written,  
\begin{equation}\label{Mom1}
\langle X^p \rangle = 
\frac{
\left. 
\left( z \frac{\partial}{\partial z} \right)^p 
\sum_{X=0}^{\infty} \frac{(X+X_A)!(X+X_B)!}{X!(X+N_f)!} z^X 
\right|_{z=1}
}
{
\left.
\sum_{X=0}^{\infty} \frac{(X+X_A)!(X+X_B)!}{X!(X+N_f)!} z^X 
\right|_{z=1}
}
\end{equation}
where the operator $(z\partial/\partial z)^p f(z) |_{z=1}$
represents  applying 
$z \times \partial /\partial z$ to $f(z)$ $p$ times, and then
evaluating the result at $z=1$. 
The denominator of (\ref{Mom1}) is simply $1/K$. 
Equation (\ref{Mom1}) differs slightly from conventional moment generating
functions \cite{Stirzaker}, in that the factor of $z$ before  
$\partial /\partial z$   
ensures that repeated application of $(z\partial /\partial z)$ yields
the moments, not the ``factorial moments'' \cite{Stirzaker} that would
be obtained by repeatedly applying $(\partial/\partial z)$.    
The hypergeometric function is defined for $|z|<1$ by \cite{Arfken}, 
\begin{equation}\label{F1}
_2F_1(a+1,b+1,c+1,z) = \frac{c!}{a!b!} 
\sum_{n=0}^{n=\infty} 
\frac{(n+a)!(n+b)!}{n!(n+c)!} z^n 
\end{equation}
provided $c\neq 0$, $-1$, $-2$, ... . 
It also has an integral representation \cite{Arfken}, 
\begin{equation}\label{F2}
_2F_1(a+1,b+1,c+1,z) =
\frac{c!}{b!(c-b-1)!} 
\int_0^1 t^b (1-t)^{c-b-1} (1-tz)^{-a-1} dt 
\end{equation}
that is valid for $|z|<1$ and $z=1$ provided Re$(c+1)>$ Re$(b+1)>0$. 
This standard result (\ref{F2}) is not obviously symmetric with
respect to $a$ and $b$ as would be expected from (\ref{F1}), however 
the expected symmetry 
is recovered later in (\ref{av1}) and (\ref{sig1}) when the
calculation is complete. 
As a consequence of (\ref{F1}), (\ref{Mom1}) can be written as, 
\begin{equation}\label{2E14}
\langle X^p \rangle = 
\frac{
\left. 
\left( z \frac{\partial}{\partial z} \right)^p 
{_2F_1}(X_A+1,X_B+1,N_f+1,z)
\right|_{z=1}
}
{
\left.
{_2F_1}(X_A+1,X_B+1,N_f+1,z)
\right|_{z=1}
}
\end{equation}
with the requirements of Re$(N_f+1)>$ Re$(X_B+1)>0$, clearly
satisfied. 
Equation (\ref{2E14}) is easily evaluated. 
Firstly  use (\ref{F2}) to substitute
for ${_2F_1}(X_A+1,X_B+1,N_f+1,z)$, then take derivatives, and set 
$z=1$. 
The resulting integral can be evaluated using 
the beta function's identity \cite{Arfken},  
\begin{equation}\label{F3} 
\int_0^1 t^b (1-t)^{c-b-1} (1-t)^{-a-1} dt 
= \frac{b!(c-a-b-2)!}{(c-a-1)!} 
\end{equation}
that holds provided Re$(c+1)>\mbox{Re}(a+1) + \mbox{Re}(b+1)$ and
Re$(b+1)>0$, a requirement that will 
restrict the values of $N_{AB}$ for 
which the resulting formulae can be used. 
This is relatively straightforward because for $t \in (0,1)$ and
  $|z|\leq 1$, $(1-tz)^{-a-1}$ is continuous with
respect to both 
$t$ and $z$, and we can bring the derivative with respect to $z$
inside the integral. Then  noting that,  
\begin{equation}\label{g1}
z \frac{\partial}{\partial z} 
\left(\frac{1}{(1 -tz)}\right)^{a +1} 
= 
\frac{a+1}{(1-tz)^{a+2}} - \frac{a+1}{(1-tz)^{a+1}} 
\end{equation}
and applying $z \partial/\partial z$ to (\ref{F2}) $p$ times, we get, 
\begin{eqnarray}\label{g2}
\lefteqn {\left. \left( z \frac{\partial}{\partial z} \right)^p 
{_2F_1}(a+1,b+1,c+1,z)
\right|_{z=1} = } & \nonumber \\ 
& & \qquad \qquad \qquad \quad (a + 1) \bigg( z
\frac{\partial}{\partial z} \bigg)^{p-1}  
\bigg[ 
{_2F_1}(a+2,b+1,c+1,z) 
\nonumber \\
& & \qquad \qquad \qquad \qquad \qquad \qquad \quad \left. {} -
	  {_2F_1}(a+1,b+1,c+1,z)  
\bigg]\right|_{z=1} 
\end{eqnarray}
where the use of (\ref{g1}) can be seen by setting $p=1$. 
Equation (\ref{g2}) 
can be iterated until the right hand side is a function of
${_2F_1}(a,b,c,1)$, for various $a$'s, $b$'s, and $c$'s, and can be
evaluated using (\ref{F3}). For $\langle X \rangle$ this gives the
average number of items missed as,  
\begin{equation}\label{av1}
\langle X \rangle = 
\frac{(N_A-N_{AB}+1)(N_B-N_{AB}+1)}{\left( N_{AB}-2\right)} 
\mbox{ with }N_{AB}>2 
\end{equation}
where $X_A$, $X_B$, and $N_f$ have been written in terms of $N_A$,
$N_B$, and $N_{AB}$, and $N_{AB}>2$ arises from the requirement on
$a$, $b$, and $c$, that allows (\ref{F3}) to be used. Similarly the
standard deviation $\sigma$ is found from, 
\begin{equation}\label{sig1}
\sigma^2 = 
\frac{(N_A-N_{AB}+1)(N_B-N_{AB}+1)(N_A-1)(N_B-1)} 
{\left(N_{AB}-2\right)^2\left(N_{AB}-3\right)} 
\mbox{ with }N_{AB}>3
\end{equation} 
Higher moments are also easily calculated and 
expressions for the skewness and kurtosis are given in the 
online supplementary material. 
Equations (\ref{av1}) and (\ref{sig1}) are exact under the assumptions
for which the prior $P(N|N_A,N_B)$ in (\ref{neweq}) 
does not depend on $N$.  
The constraints on the minimum
value of $N_{AB}$ for which the expressions hold is a mathematical
requirement, and appears to be a requirement for the series to
converge. As discussed later, this requirement on $N_{AB}$ can be
overcome with a suitably convergent prior distribution $P(N)$. Because
both $N_A$ and $N_B$ are greater than or equal to $N_{AB}$, then
$N_{AB}>2$ will require $N_A > 2$ and $N_B > 2$ also.


\subsection{Comparison with Chapman's estimate} 

Previous approaches have approximated these same average and standard
deviation by a combination of conjecture and estimations for the
precision and bias \cite{1951,Seber1970a,Wittes,Seber}. It has been
observed \cite{Pelayo} that previous (approximate) estimates can be
inaccurate for combinations of $N_A$, $N_B$, and $N_{AB}$ that cause
the hypergeometric distribution  to have a `long tail', for example if
$N_A \gg N_B$. These remarks can now be 
clarified. \citename{1951}'s (1951) 
estimation gives $\langle N \rangle \approx
\frac{(N_A+1)(N_B+1)}{(N_{AB}+1)}-1$, and $\langle X \rangle = \langle
N \rangle - N_f$, as, 
\begin{equation}\label{Cest}
\langle X \rangle \approx 
\frac{(N_A-N_{AB})(N_B-N_{AB})}{(N_{AB}+1)} 
\end{equation}
Comparing this
with (\ref{av1}) (for example by subtracting (\ref{Cest}) from
(\ref{av1})), we  can see that:  
\begin{enumerate}
\item it is always less than (\ref{av1}),
\item that this is more pronounced when either or both of $(N_A-N_{AB})$ or  
  $(N_B-N_{AB})$ are large, or when $N_{AB}$ is small, but that
  conversely,   
\item provided neither $N_A$ nor $N_B$ equals $N_{AB}$, it will give
  the same (unbiased) estimate if $N_{AB}$ is sufficiently large
  compared with both $(N_A-N_{AB})$ and $(N_B-N_{AB})$.   
\end{enumerate}
Similar remarks apply to the widely used estimate for the variance
\cite{Seber1970a}, that has 
\begin{equation}
\sigma^2 \approx
\frac{(N_A+1)(N_B+1)(N_A-N_{AB})(N_B-N_{AB})}
{(N_{AB}+1)^2(N_{AB}+2)}
\end{equation} 
and is
unbiased for $N_{AB} \gg 1$, but accuracy requires an increasingly
large $N_{AB}$ if either $(N_A-N_{AB})$ or $(N_B-N_{AB})$ are small,
and in practice it can be inaccurate. 

Seber \citeyear{Seber1970a,Seber} has remarked that \citename{1951}'s
calculations are equivalent to approximating (\ref{neweq}) with a
Poisson distribution.
Appendix B 
finds this requires both 
$0 \neq (N_A-N_{AB})/N_{AB} \ll 1$ and  $0\neq (N_B - N_{AB})/N_{AB}
\ll 1$, (and implicitly that $N_{AB} \gg 1$). 
When this is true, the mean of the approximating Poisson
distribution coincides with the maximum of (\ref{Papp}) 
with $\langle X \rangle =(N_A-N_{AB})(N_B-N_{AB})/N_f$, and approximates
both (\ref{av1}) and (\ref{Cest}) (for this limit). 
Similarly for the variance. 
In
contrast (\ref{av1}) and (\ref{sig1}) result from exactly calculating
the moments of (\ref{Papp}).  
As noted in Appendix B,  
this Poisson approximation 
generalises to the situation studied by \cite{Pelayo}, in which there
are $n$ searches instead of only two. 


\subsection{The moments of (\ref{fulleq})}\label{3.3}

When all items are searched for by both A and B, the probability
distribution for the number of items searched for is given by
(\ref{fulleq}). 
For the common choice of prior with $P(N)$ constant, 
Appendix C 
shows how the moments of (\ref{fulleq}) can be
closely approximated using the moments of (\ref{neweq}), and
calculates rigorous maximum bounds for the error in the approximation.  
When $N_f \gg 1$ the error will be small and a good approximation is 
given by,   
\begin{equation}\label{av2}
\langle X \rangle =   
\frac{(N_A-N_{AB}+1)(N_B-N_{AB}+1)}{N_{AB}} 
\mbox{ for }N_{AB}>0
\end{equation}
with an error that is less than $\pm \langle X \rangle/(N_f+1)$.  
Unfortunately $\sigma^2=\langle X^2 \rangle - \langle X \rangle^2$ can
be arbitrarily small, but the approximation for $\sigma^2$ of,  
\begin{equation}\label{sig2}
\sigma^2 =
\frac{(N_A-N_{AB}+1)(N_B-N_{AB}+1)(N_A+1)(N_B+1)} 
{N_{AB}^2\left(N_{AB}-1\right)} 
\mbox{ with }N_{AB}> 1
\end{equation}
has a maximum error that is of order $\langle X^2 \rangle/N_f$. 
Consequently unless $\langle X^2 \rangle/N_f \ll 1$, (\ref{sig2}) is not
guaranteed to be a good approximation for $\sigma^2$.  
Often there will be a prior reason to expect that $N \gg 1$. 
For these cases an alternative approach is to assume the almost
constant prior of,  
\begin{equation}\label{kappa}
P(N)=\kappa \frac{(N+1)}{(N+2)} 
\end{equation}
with $\kappa$ constant, that may be written as $P(N) = \kappa (1
-1/(N+2))$, and monotonically increases from $P(0)=\kappa/2$ to
$P(\infty) = \kappa$.  
This prior gives a small bias against low values of $N$ but is
approximately constant for larger values of $N$.  
For example, $P(N)$ varies by less than ten percent between $N=8$ 
and $N=\infty$. 
For this prior (\ref{fulleq}) becomes, 
\begin{equation}\label{hash} 
P(N|N_A,N_B,N_{AB}) 
= 
\frac{ (N-N_A)!(N-N_B)! }{
(N+2)!(N-N_f)!}
\kappa
\end{equation}
Remembering that $N_f=N_A+N_B-N_{AB}$, then rewriting (\ref{hash}) in
terms of $(N+2)$, $(N_A+2)$, $(N_B+2)$, and $(N_{AB}+2)$, it will be
clear 
that the change of variables that replaces: $(N+2)$ with $N$, 
$(N_A+2)$ with $N_A$, $(N_B+2)$ with $N_B$, $(N_{AB}+2)$ with
$N_{AB}$, 
makes (\ref{hash}) the same form as (\ref{neweq}).
The condition that $N=N_f$ may be written as 
$(N+2)=(N_A+2)+(N_B+2)-(N_{AB}+2)$, so after the change of variables
the lower limit $N=N_f$ on sums for the moments remains the same.  
The upper limit of $N=\infty$ is clearly also unchanged. 
Consequently the exact moments of (\ref{hash}) can be found by replacing
$N_A$ with $N_A+2$, $N_B$ with $N_B+2$, and $N_{AB}$ with $N_{AB}+2$,
in the exactly calculated moments of (\ref{neweq}), with for example
(\ref{av1}) and (\ref{sig1}) becoming (\ref{av2}) and (\ref{sig2}). 
(An alternative presentation of these remarks can be found in Appendix
C.) 
With the prior (\ref{kappa}), 
(\ref{av2}) and (\ref{sig2}) are exact moments of (\ref{fulleq}), and 
the error bounds now provide a bound on the maximum possible difference
between estimates calculated with this, and with a  flat prior. 
For those cases when it is reasonable to assume this prior, we think
it is preferable to explicitly use it along with the exact estimates
(\ref{av2}) and (\ref{sig2}), in  preference to assuming a constant
prior and treating (\ref{av2}) and (\ref{sig2}) as approximations.


Both (\ref{av2}) and (\ref{sig2}) are more similar to the Chapman and
Lincoln-Petersen estimates than (\ref{av1}) and (\ref{sig1}). 
This is despite them
being approximations to the moments of (\ref{fulleq}), not
(\ref{neweq}), that Chapman's calculation is intended to approximate. 
This might help explain why the discrepancy between Chapman's 
estimate and (\ref{av1}) is generally overlooked. 
For many cases of interest the number of items found ($N_f$) is large,
with $N_f \gg1$, and for these cases (\ref{av2}) 
provides an accurate estimate for $\langle X \rangle$. 
Next we consider some examples.

\subsection{Examples}

When A and B each search for a number of items that is   
predetermined in advance of their search, then (\ref{av1}) and (\ref{sig1}) 
provide simple estimates for the maximum number of items that could 
be found by a search for all items, and the precision of the estimate.
They are exact moments of (\ref{neweq}). 
When all items are searched for, provided the number of items found
($N_f$) is much greater than one, then a very good estimate can be made
using (\ref{av2}), and if the prior $P(N)=\kappa (N+1)/(N+2)$ is
assumed then (\ref{av2}) and (\ref{sig2}) are exact moments of
(\ref{fulleq}).    
Both pairs of estimates can  
give substantially different estimates to those
of Chapman (\ref{Cest}) and Lincoln-Petersen (\ref{X1}). 
For example, \citeasnoun**{Chao} propose a method to combine multiple
intersections of lists and the Lincoln-Petersen or Chapman estimator,
with the intention of improving the accuracy of epidemiological 
estimates.  
The number of items in common between lists is not predetermined, and
is anywhere between zero and every item on the shortest list.  
Their proposed method is illustrated in Section 4 of
\citeasnoun**{Chao}, and
the estimates calculated by the method are given on the top of page
968, where they are calculated from the numbers in their Table 5b
using the 
Chapman and also the Lincoln-Petersen estimate.  
The results of their
calculations are reported in Table 6 on page 968 of their paper, and
repeated in part in Table \ref{table1}.   
The total number of items ($N_f$) is much larger than one in all
cases, and consequently an accurate estimate is given by (\ref{av2}).  
An immediate concern is that the Chapman and Lincoln-Petersen estimates are
estimators for the moments of (\ref{neweq}), that arise from a search
procedure for a predetermined number of items, and should not be used.  
It is a fortunate coincidence that the moments of (\ref{fulleq}) are
closer to the Chapman and Lincoln-Petersen estimates than are the
exact moments of (\ref{neweq}) that they are intended to approximate. 
They are also estimates 
for the most probable population size, and not the
expectation of the population size, which can be much larger. 
For the cases in Table 5b of \citeasnoun**{Chao} where (\ref{av2}) and  
(\ref{sig2}) are defined, we find the revised estimates given in Table
\ref{table1}.   
\begin{table}
\begin{center}
\begin{tabular}{|l|c|c|c|c|c|c|c|c|}
\hline
{ } & $N_A$ & $N_B$ & $N_{AB}$ 
& $\langle N \rangle$ & $\sigma$ 
& $\langle N \rangle_{C}$ & $\langle N \rangle_{LP}$ & $\sigma_{S}$ \\
\hline
Male & 323 & 101 & 3 
& 11014 
& 7638 
& 8261 & 10874 & 3599 \\
\hline
Female & 21 & 19 & 1
& 438 
& undefined 
& 219 & 399 & 115 \\
\hline 
Combined & 344 & 120 & 4
& 10434 
& 5890 
& 8348 & 10320 & 3067 \\
\hline
\end{tabular}
\end{center}
\caption[]{\label{table1} Estimates for $\langle N \rangle = N_f +
  \langle X \rangle$ and $\sigma$ are calculated using (\ref{av2}),
  (\ref{sig2}), and the numbers in Table 5b of
  \citeasnoun**{Chao}, that are reproduced above as $N_A$, $N_B$, and
  $N_{AB}$.  
The estimates from Table 6 of \citeasnoun**{Chao},
that use the Chapman ($\langle N \rangle_{C}$) and Lincoln-Petersen
  ($\langle N \rangle_{LP}$) estimates for $\langle N
\rangle$, and Seber's estimate for the variance ($\sigma_S$), are also
  included. 
Our estimates, where they are defined, are substantially different to
  the quoted estimates {\cite{Chao}} that use the
  Lincoln-Petersen (\ref{X1}) and Chapman (\ref{Cest})
  estimates. 
} 
\end{table}
Also included are the estimates from Table 6 of \citeasnoun**{Chao},
and Seber's estimate for the variance. 
Our estimates are substantially different, and in some cases $N_{AB}$
is too small to allow them to be used. 
It is unusual, but not unreasonable, to find distribution
functions without a well-defined mean or standard deviation. 
Without a suitable prior distribution the 
female list for the ``shared population'' of \citeasnoun**{Chao} will fall  
into this category.  
For such cases it is necessary to (explicitly) use a suitable prior if
estimates are to be correctly made.   

Smaller deviations from the usual Lincoln-Petersen and Chapman
estimates are expected when $N_{AB}$ is sufficiently large compared to
$N_A$ and $N_B$. 
For example, in a
recent review by \citeasnoun**{May}, there were $177$ relevant papers
found by author A, $265$ papers found by author B, and $171$ of these
papers found by both authors (K.E. May, private communication). Using
(\ref{av2}) and (\ref{sig2}), we find $\langle X \rangle \simeq 3.9$
and $\sigma = 2.5$. 
Therefore whereas $271$ papers were found, our 
estimate gives between $1$ and $6$ missed papers. 
Putting it
another way, the estimate is that between $97.6\%$ and $99.5\%$ of the
papers searched for from within the total sample of just over $8$
thousand papers were found. 
The standard estimates
\cite{1951,Seber1970a} give $\langle X \rangle = 3.3$ and $\sigma =
2.3$, and are somewhat 
smaller despite the reasonably large value of $N_{AB}=171$. 
Another literature search example \cite{Spoor} found  $N_A=150$,
$N_B=123$, and $N_{AB}=115$, for which (\ref{av2}) and (\ref{sig2})
give $\langle X \rangle = 2.8$ and $\sigma = 2.0$. These compare with
the standard estimates \cite{1951,Seber1970a} that give, $\langle X
\rangle = 2.4$ and $\sigma = 1.8$.

\subsection{Limitations of the model}

Underlying the calculation is the assumption that all items are
equally likely to be found. 
Clearly there will be cases where some items are more difficult to
find. 
However even in those cases, some (lower bound) estimate for the
number of items missed is better than no estimate at all. 
The method will fail most dramatically  if there is a sub-population
that is much more difficult to find; it is possible that both
searchers could miss all or most of that sub-population, and will
overestimate the accuracy of their search. 
These limitations should be considered before applying these
estimates, and when reporting them. 
If there is a (prior) reason to think the assumptions are
inappropriate, one way that modified assumptions can be included is
through different  priors for $P(N_A|N)$ and $P(N_B|N)$ as was
discussed in Section \ref{S2.1}.
In general this will give distribution functions that are 
most easily calculated numerically. 
 

\section{Bayesian corrections and other search procedures}\label{sbayes2} 

An advantage of the Bayesian approach is that the assumptions are
explicit at the outset and the resulting answers are exact, with no
additional free parameters. Before concluding we consider two easily
evaluated examples that illustrate how different prior assumptions and 
different search procedures affect the estimates. 

\subsection{One partial and one comprehensive search}

Firstly imagine a situation where one author (\emph{e.g.} A) searches
for a fixed number of papers so that $P(N_A|N)$ no longer appears in
(\ref{fulleq}),  
but the other author (B)  searches for as many papers as possible with
$P(N_B|N)=1/(N+1)$, with no prior knowledge of the number of papers
searched for other than it being finite ($P(N)$ constant). 
For this case (\ref{neweq}) is modified by the factor $1/N!$ becoming 
$1/(N+1)!$. 
In Section \ref{3.3} it was explained how a suitable change of
variables could transform (\ref{hash}) into the same form as
(\ref{neweq}), allowing the moments of (\ref{hash}) to be calculated
from those of (\ref{neweq}) by a simple change of variables.  
The same is true here, the change of variables that replaces: $(N+1)$
with $N$, $(N_A+1)$ with $N_A$, $(N_B+1)$ with $N_B$, $(N_{AB}+1)$
with $N_{AB}$, leads to the same form of $P(N|N_A,N_B,N_{AB})$ as
(\ref{neweq}). 
Similarly to Section \ref{3.3}, because the equation
$N=N_f=N_A+N_B-N_{AB}$ may be written as
$(N+1)=(N_A+1)+(N_B+1)-(N_{AB}+1)$, the lower limit on the range of
summation for the moments remains unchanged by the change of
variables, as does the $N=\infty$ upper limit.  
Consequently the exact moments can be found by replacing $N_A$ by
$N_A+1$, $N_B$ by $N_B+1$, $N_{AB}$ by $N_{AB}+1$, in (\ref{av1}) and
(\ref{sig1}), giving,   
\begin{equation}
\langle X \rangle = 
\frac{(N_A-N_{AB}+1)(N_B-N_{AB}+1)}{\left( N_{AB}-1\right)} 
\mbox{ with }N_{AB}>1 
\end{equation}
and, 
\begin{equation}
\sigma^2 = 
\frac{(N_A-N_{AB}+1)(N_B-N_{AB}+1)(N_A)(N_B)} 
{\left(N_{AB}-1\right)^2\left(N_{AB}-2\right)}
\mbox{ with }N_{AB}>2
\end{equation}


Interestingly, for this search procedure the standard
capture-recapture estimate 
conjectured by \citename{1951} of $\langle N \rangle \approx
\frac{(N_A+1)(N_B+1)}{(N_{AB}+1)}-1$, approximates the ``most
probable'' value of $N$, where $P(N|N_A,N_B,N_{AB})$ is a maximum. The
maximum can be approximated by setting
$P(N|N_A,N_B,N_{AB})=P(N-1|N_A,N_B,N_{AB})$ and solving for $N$
\cite{1951,Pelayo}. For the stated prior assumptions this gives,  
\begin{equation}
\frac{(N-N_A)!(N-N_B)!}{(N+1)!(N-N_A-N_B+N_{AB})!} 
= 
\frac{(N-N_A-1)!(N-N_B-1)!}{N!(N-N_A-N_B+N_{AB}-1)!} 
\end{equation}
whose solution for $N$ is exactly Chapman's 
conjectured estimate. (Strictly this estimate
is only an approximation to the most probable value of $N$: a more
precise value can be found using Stirling's approximation for the
factorials and differentiating with respect to $N$ to find the maximum
of $P(N|N_A,N_B,N_{AB})$.) 


\subsection{The influence of a proper prior}

To illustrate the effect of $P(N)$, consider the normalisable prior
$P(N)=\kappa (N+1)/(N+2)(N+3)(N+4) \sim \kappa/N^2$, with $\kappa$
constant, and let both A and B search 
for as many items as possible with $P(N_A|N)=P(N_B|N)=1/(N+1)$. For
this example (\ref{neweq}) is modified by $1/N!$ becoming
$1/(N+4)!$. 
Following a similar change of variables as discussed above and in
Section \ref{3.3}, but now with: $(N+4)$ replaced by $N$, $(N_A+4)$
with $N_A$, $(N_B+4)$ with $N_B$, $(N_{AB}+4)$ with $N_{AB}$, then
$P(N|N_A,N_B,N_{AB})$ becomes the same form as in (\ref{neweq}).  
Consequently modified estimates can be found by substituting $N_A$
with $N_A+4$,  $N_B$ with $N_B+4$, and $N_{AB}$ with $N_{AB}+4$, in
(\ref{av1}) and (\ref{sig1}), leading to a reduced estimate for
$\langle X \rangle$.  

Notice that for this latter example the requirement that $N_{AB}>3$ in
(\ref{sig1}) becomes (with $N_{AB}$ replaced by $N_{AB}+4$), $N_{AB}
> -1$, and the estimates hold for all $N_A$, $N_B$,
and $N_{AB}$. The conclusion is that  whereas (\ref{av1}) and
(\ref{sig1}) can only be used when $N_{AB}$, $N_A$, and $N_B$ are
sufficiently large ($>3$), when all items are searched for (resulting
in the extra factor of $1/(N+1)^2$ in $P(N|N_A,N_B,N_{AB})$), the 
equations apply for a greater range of values. In fact 
unless $N_{AB}$ is sufficiently large, 
then estimates can {\sl only} be calculated with a sufficiently
convergent ({\sl i.e.} realistic) prior for a given search strategy
(such as searching for a fixed number of items, or for all the items).  
In summary, it is important to ensure that the assumptions upon which
any given estimate depends are consistent with the problem being
studied. 


\section{Conclusions}\label{conc}

The original purpose of this calculation was to consider
two authors A and B searching a finite set of papers for those to
include in a literature survey, and to use the number of papers found
by authors A ($N_A$) and B $(N_B$), 
along with the number found by {\sl both}
authors ($N_{AB}$), to estimate how accurate the search was. Bayes'
theorem is used to rigorously formulate 
this ``mark-recapture'' problem for two different
search procedures. 
The first procedure corresponds to A and B searching for all of the
items, the second corresponds to A and B each searching for a
predetermined number of items, before comparing their
results to allow an estimate for $N$.   
For the latter case, exact
calculations lead to simple formulae for the average number of items 
missed from the search (\ref{av1}), and the standard deviation
(\ref{sig1}). 
The skewness and kurtosis of the probability
distribution are given within the appendices in the online 
supplementary information, 
and higher moments may be
calculated in a similar way. 

Equations (\ref{av1}) and (\ref{sig1}) are exact moments of the
widely-studied 
probability distribution (\ref{neweq}) from Chapman's 1951 paper,
which is shown here to result from a procedure in
which A and B each search for a {\sl predetermined} number of items.
Previous estimates using this distribution 
have been derived using a combination of conjecture and approximations. 
Chapman's 
conjectured estimate is found (under suitable assumptions) to
be an approximation to the most probable value of $N$. This provides a
good approximation to (\ref{av1}) if $N$ is large and both
searchers individually find the majority of the items searched for, but
is increasingly bad if either searcher finds substantially more (or
fewer) items than their partner, which can often be the case.  

For many cases such as the literature search application, all items
are searched for by both A and B, which leads to a modified
probability distribution (\ref{fulleq}). 
If a constant prior is assumed then the moments of (\ref{fulleq}) can be 
closely approximated 
provided the
number of items found ($N_f$) is much greater than one, which will
very often be the case.  
When this is the case, an excellent approximation for the number of
items missed is given by (\ref{av2}).  
Alternately if there is a prior reason to think $N \gg 1$, then it is
reasonable to use the almost constant prior 
$P(N)=\kappa (N+1)/(N+2)$,
and the calculation for the estimates of (\ref{av2}) and (\ref{sig2})
becomes  exact.  
For estimates arising from this search procedure, there is a smaller
difference 
between them and Chapman's estimate (which we have shown here does not
apply, and in principle should not be used), 
but it can still be 
substantial. We recommend using the improved estimates given by 
(\ref{av1}), (\ref{sig1}), (\ref{av2}), and (\ref{sig2}), as is
appropriate to the search procedure.

The formulae apply to an enormously wide variety of problems
with two independent searches in which the number of items found by
searcher A ($N_A$), searcher B ($N_B$), and the number found by both
($N_{AB}$), can be determined. By ``independent'', we mean that A
finding an item does not affect the probability of B finding it
(\emph{e.g.} for mark-and-recapture, animals do not become ``shy'' or
``tame'' after handling). Finally we caution against an assumption
used in the calculation -- that all objects searched for are equally
likely to be found. This will fail if there is a sub-population that
is much more difficult to find, for which case both searchers will
appear to have found the majority of items and will over-estimate the
accuracy of their search. 
These issues are beyond the intended scope of this paper. 
Nonetheless even when the assumption is only
approximately true (often the assumption will be good), 
these improved estimates (\ref{av1}), (\ref{sig1}), (\ref{av2}), and
(\ref{sig2}) will hopefully 
provide a valuable  standard 
tool for literature searches and more generally.


\section*{Acknowledgements}
Thanks to Dr Katie Webster (previously Dr Katie May) for recording and
supplying the numbers from 
the recent literature search described in \citeasnoun{May}, for
emphasising the potential use of this technique in literature
searches, and for numerous helpful discussions. Thanks to Professor
Walter Zucchini for supplying a copy of \citeasnoun{Zucchini}, and
Martin O'Brien for helpful discussions and comments. 
Thanks also to the editor of The American Statistician, and the associate
editor in particular, for numerous  
helpful comments and suggestions.

\appendix

\section{The moments}\label{moments}


Here we briefly present our original derivation of the moments of 
(10),  
that uses simpler mathematical concepts, but is less conventional and
systematic than the
generating function approach presented in the main text. 
Repeating (10) 
here for convenience, with, 
\begin{equation}\label{e30} 
\begin{array}{c}
P(X|X_A,X_B,N_{AB}) =  
\frac{(X+X_A)!(X+X_B)!}{X!(X+N_{f})!} 
K
\end{array} 
\end{equation}
and $X$ between $0$ and $\infty$. 
Next define,
\begin{equation}
\begin{array}{c}
S(X_A,X_B,N_{f}) =  
\sum_{X=0}^{\infty} \frac{(X+X_A)!(X+X_B)!}{X!(X+N_{f})!} 
\end{array} 
\end{equation}
where we note that $N_f=N_{AB}+X_A+X_B$, and also that
$K=1/S(X_A,X_B,N_{f})$.  
The aim is to express the moments in terms of the function
$S(X_A,X_B,N_f)$, evaluate $S(X_A,X_B,N_f)$ using an identity due to
Gauss, then combine the results to obtain explicit expressions for the
moments in terms of $X_A$, $X_B$, and $N_f$.

Starting with $\langle X \rangle$, notice that, 
\begin{equation}
\begin{array}{ll}
\sum_{X=0}^{\infty} X 
\frac{
\left(X+X_A\right)! \left(X+X_B\right)! }
{
X! \left(X+N_f\right)! 
}
&= \sum_{X=1}^{\infty}  \frac{X}{X!} 
\frac{
\left(X-1+X_A+1\right)! \left(X-1+X_B+1\right)! }
{
\left(X-1+N_f+1\right)! 
}
\\
&= \sum_{X=1}^{\infty}  \frac{1}{(X-1)!} 
\frac{
\left((X-1)+(X_A+1)\right)! \left((X-1)+(X_B+1)\right)! }
{
\left((X-1)+(N_f+1)\right)! 
}
\\
&= \sum_{X=0}^{\infty}   
\frac{
\left(X+X_A+1\right)! \left(X+X_B+1\right)! }
{
X! \left(X+N_f+1\right)! 
}
\\
&=S(X_A+1, X_B+1,N_f+1) 
\end{array}
\end{equation}
Hence, 
\begin{equation}\label{avX1}
\langle X \rangle = \frac{S(X_A+1, X_B+1,N_f+1)}
{S(X_A, X_B,N_f)}
\end{equation}
Similarly for $\langle X^2 \rangle$, 
\begin{equation}
\begin{array}{ll}
\sum_{X=0}^{\infty} X^2 
\frac{
\left(X+X_A\right)! \left(X+X_B\right)! }
{
X! \left(X+N_f\right)! 
}
&= \sum_{X=1}^{\infty}  \frac{X}{X!} \left( X-1+1\right)
\frac{
\left(X-1+X_A+1\right)! \left(X-1+X_B+1\right)! }
{
\left(X-1+N_f+1\right)! 
}
\\
&= \sum_{X=1}^{\infty}  \frac{\left( (X-1)+1\right)}{(X-1)!} 
\frac{
\left((X-1)+(X_A+1)\right)! \left((X-1)+(X_B+1)\right)! }
{
\left((X-1)+(N_f+1)\right)! 
}
\\
&= \sum_{X=0}^{\infty}  \left( X+1\right) 
\frac{
\left(X+X_A+1\right)! \left(X+X_B+1\right)! }
{
X! \left(X+N_f+1\right)! 
}
\end{array}
\end{equation}
Repeating the same trick to remove the factor of $X$ then gives, 
\begin{equation}\label{sigX1}
\langle X^2 \rangle = 
\frac{S(X_A+2, X_B+2,N_f+2)}
{S(X_A, X_B,N_f)}
+
\frac{S(X_A+1, X_B+1,N_f+1)}
{S(X_A, X_B,N_f)}
\end{equation}
Similarly but with more algebra for the higher order moments, e.g. 
\begin{equation}\label{skewX1}
\langle X^3 \rangle = 
\frac{S(X_A+3, X_B+3,N_f+3)}
{S(X_A, X_B,N_f)}
+
3\frac{S(X_A+2, X_B+2,N_f+2)}
{S(X_A, X_B,N_f)}
+
\frac{S(X_A+1, X_B+1,N_f+1)}
{S(X_A, X_B,N_f)}
\end{equation}
and,
\begin{equation}\label{kurtX1}
\begin{array}{ll}
\langle X^4 \rangle &= 
\frac{S(X_A+4, X_B+4,N_f+4)}
{S(X_A, X_B,N_f)}
+
6\frac{S(X_A+3, X_B+3,N_f+3)}
{S(X_A, X_B,N_f)}
\\
&+
7\frac{S(X_A+2, X_B+2,N_f+2)}
{S(X_A, X_B,N_f)}
+
\frac{S(X_A+1, X_B+1,N_f+1)}
{S(X_A, X_B,N_f)}
\end{array}
\end{equation}

To evaluate $S(X_A, X_B,N_f)$, we firstly note that the hypergeometric
function has for $|z|<1$ and $c\neq 0$, $-1$, $-2$, ... 
(Arfken 1985), 
\begin{equation}\label{B8}
_2F_1\left( a+1, b+1, c+1, z\right) 
= \frac{c!}{a!b!} \sum_{n=0}^{\infty} 
\frac{(a+n)!(b+n)!}{(c+n)!} \frac{z^n}{n!} 
\end{equation} 
For $z=1$ an identity due to Gauss gives (Arfken 1985), 
\begin{equation}\label{B9}
_2F_1\left( a+1, b+1, c+1, 1\right) 
= \frac{\Gamma(c+1) \Gamma(c-a-b-1)}{\Gamma(c-a)\Gamma(c-b)} 
\mbox{ , Re}(c) > \mbox{ Re}(a+b) + 1
\end{equation}
with $c\neq 0$, $-1$, $-2$, ... , as above.  
Equations (\ref{B8}) and (\ref{B9}) 
may be combined to give (for $z=1$), 
\begin{equation}
\sum_{n=0}^{\infty}
\frac{(a+n)!(b+n)!}{(c+n)!} \frac{1}{n!} 
= 
\frac{a!b!}{c!} 
\frac{\Gamma(c+1) \Gamma(c-a-b-1)}{\Gamma(c-a)\Gamma(c-b)} 
\mbox{ , Re}(c) > \mbox{ Re}(a+b) + 1
\end{equation}
Therefore with the replacements of $n=X$, $c=N_f$, $a=X_A$, and
$b=X_B$ (so that $c=a+b+N_{AB}>(a+b)+1$ for $N_{AB} > 1$), we get, 
\begin{equation}\label{appS}
\begin{array}{ll}
S(X_A,X_B,N_f) &=
\sum_{X=0}^{\infty} \frac{(X+X_A)!(X+X_B)!}{X!(X+N_f)!}
\\
&= X_A! X_B! \frac{ (N_f-X_A-X_B -2 )!}
{
(N_f-X_A-1)! (N_f-X_B-1)! 
}
\mbox{ , } N_{AB} > 1 
\end{array}
\end{equation}
Hence substituting into (\ref{avX1}) 
gives, 
\begin{equation}\label{avX}
\langle X \rangle = \frac{(X_A +1)(X_B+1)}{(N_f-X_A-X_B-2) }
= \frac{(X_A +1)(X_B+1)}{(N_{AB}-2) } \mbox{ , } N_{AB}>2
\end{equation}
where the inequality follows from the requirement that 
$N_f+1 > (X_A+1)+(X_B+1)+1$ with $N_f=N_{AB}+X_A+X_B$. 
Similarly, 
\begin{eqnarray}\label{avXsq}
\langle X^2 \rangle &=& 
\frac{(X_A+1)(X_A+2)(X_B+1)(X_B+2)}{(N_{AB}-2)(N_{AB}-3)} 
\nonumber \\
&&{}+ 
\frac{(X_A+1)(X_B+1)}{(N_{AB}-2)} \mbox{  , with  } N_{AB}>3
\end{eqnarray}
\begin{eqnarray}\label{skewX2}
\langle X^3 \rangle &=& 
\frac{(X_A+1)(X_A+2)(X_A+3)(X_B+1)(X_B+2)(X_B+3)}
{(N_{AB}-2)(N_{AB}-3)(N_{AB}-4)}
\nonumber \\
&&{}+
3\frac{(X_A+1)(X_A+2)(X_B+1)(X_B+2)}
{(N_{AB}-2)(N_{AB}-3)}
\nonumber \\
&&{}+
\frac{(X_A+1)(X_B+1)}
{(N_{AB}-2)}   \mbox{  , with  } N_{AB}>4
\end{eqnarray}
\begin{eqnarray}\label{kurtX2}
\langle X^4 \rangle &=& 
\frac{(X_A+1)(X_A+2)(X_A+3)(X_A+4)(X_B+1)(X_B+2)(X_B+3)(X_B+4)}
{(N_{AB}-2)(N_{AB}-3)(N_{AB}-4)(N_{AB}-5)}
\nonumber \\
&&{}+
6\frac{(X_A+1)(X_A+2)(X_A+3)(X_B+1)(X_B+2)(X_B+3)}
{(N_{AB}-2)(N_{AB}-3)(N_{AB}-4)}
\nonumber \\
&&{}+
7\frac{(X_A+1)(X_A+2)(X_B+1)(X_B+2)}
{(N_{AB}-2)(N_{AB}-3)}
\nonumber \\
&&{}+
\frac{(X_A+1)(X_B+1)}
{(N_{AB}-2)}
\mbox{  , with  } N_{AB}>5
\end{eqnarray}
These may be used to calculate various statistical quantities. 
The standard deviation 
$\sigma = \sqrt{ \langle X^2 \rangle -\langle X \rangle^2 }$, which using 
(\ref{avX}) and  (\ref{avXsq}), 
simplifies to give, 
\begin{equation}\label{sigX}
\sigma 
= \sqrt{
\frac{(X_A+1)(X_B+1)(N_{AB}+X_A-1)(N_{AB}+X_B-1)}
{(N_{AB}-2)^2(N_{AB}-3)}
}   \mbox{ , } N_{AB}>3
\end{equation}
The skewness $\gamma = 
\langle 
\left( X - \langle X \rangle \right)^3 \rangle 
/ \langle X^2 \rangle^{3/2}
$, which expands to give, 
\begin{equation}
\gamma = \frac{
\langle X^3 \rangle 
- 3 \langle X^2 \rangle \langle X \rangle 
+ 2 \langle X \rangle^3
}
{
\langle X^2 \rangle^{3/2}
}
\end{equation}
and may be evaluated using (\ref{avX}) to (\ref{skewX2}). 
The kurtosis is given by $\kappa = 
\langle 
\left( X - \langle X \rangle \right)^4 \rangle 
/ \langle X^2 \rangle^{2}
$, which expands to give, 
\begin{equation}\label{kurt}
\kappa = 
\frac{
\langle X^4 \rangle - 4 \langle X^3 \rangle \langle X \rangle 
+6 \langle X^2 \rangle \langle X \rangle^2 - 3 \langle X \rangle^4
}
{
\langle X^2 \rangle^{2}
}
\end{equation}
and may be evaluated using  (\ref{avX}) to (\ref{kurtX2}). 
Replacing $X=N-N_f$, $X_A=N_A-N_{AB}$, and $X_B = N_B - N_{AB}$
in (\ref{avX}) and (\ref{sigX}), 
gives 
(19) 
and 
(20) 
of the main text.

\section{Poisson approximation}\label{app2}

Starting from 
(11) in the main text, 
use the approach of 
Chapman (1951) and Garc\'ia-Pelayo (2006) 
to find $X$ for which $P(X|X_A,X_B,N_{AB})$ is
maximum, from $P(X^*|X_A,X_B,N_{AB})=P(X^*-1|X_A,X_B,N_{AB})$.  
This gives $X^*=X_AX_B/N_f=(N_A-N_{AB})(N_B-N_{AB})/N_f$. 
When both $X_A/N_{AB} \ll 1$ and $X_B/N_{AB} \ll 1$, then both $X^*
\ll X_A$ and $X^* \ll X_B$, and because $N_f=N_{AB}+X_A+X_B$ is larger
than either $X_A$ or $X_B$ then $X^* \ll N_f$ also.

Next note that $(X+X_A)!\equiv X_A! X_A^X \exp\{\sum_{y=1}^X
\log(1+y/X_A)\}$, as may be seen from expanding $(X+X_A)!$, 
\begin{equation}
\begin{array}{ll}
(X+X_A)! &= (X_A + X)(X_A + X -1) ... (X_A+1)X_A!
\\
&=X_A!\exp\left\{\sum_{y=1}^X \log(y+X_A)\right\}
\end{array}
\end{equation}
where the last line repeatedly used $AB=\exp(\log(AB))=\exp(\log(A)+\log(B))$. 
Then write,
\begin{equation}
\begin{array}{ll}
X_A!\exp\left\{\sum_{y=1}^X \log(y+X_A)\right\} 
&= X_A!\exp\left\{\sum_{y=1}^X \log(X_A(1+y/X_A))\right\}
\\
&= X_A!\exp\left\{X\log(X_A) + \sum_{y=1}^X \log(1+y/X_A)\right\}
\\
&= X_A! \exp\left\{ \log(X_A^X)\right\} \exp\left\{\sum_{y=1}^X
  \log(1+y/X_A) \right\} 
\\
&= X_A! {X_A}^{X} \exp\left\{\sum_{y=1}^X \log(1+y/X_A)\right\}
\end{array}
\end{equation}
as originally stated. 
Similarly expanding $(X+X_B)!$ and $(X+N_f)!$, gives, 
\begin{equation}\label{exact2}
\begin{array}{ll}
P(X|X_A,X_B,N_f) &= K \frac{X_A!X_B!}{N_f!} \frac{1}{X!} 
\left( \frac{X_AX_B}{N_f} \right)^{X} 
\times
\\
&\exp\left\{ 
\sum_{i=1}^{X} \log(1+i/X_A) + \sum_{j=1}^X \log(1+j/X_B) 
- \sum_{k=1}^X \log(1+k/N_f) 
\right\} 
\end{array}
\end{equation} 

The above expression is exact, and can be used as the starting point
for a variety of approximations.  
It is composed of the product of a Poisson distribution
${X^*}^X/X!$ with $X^*=X_AX_B/N_f$, a constant term that ensures
(\ref{exact2})   
is correctly normalised, and an exponential
term whose exponent is a function of $X$.  
As $X$ becomes small relative to $X_A$, $X_B$, and $N_f$, the
exponential's exponent tends to zero, and (\ref{exact2}) 
asymptotes to a Poisson distribution. 
However, because $X_A<N_f$ and $X_B<N_f$, the exponential term's
exponent is a strictly increasing function of $X$. 
Consequently a good approximation to (\ref{exact2}) 
by a
Poisson distribution is only ever possible over a limited range of $X$. 
An approximation with a Poisson
distribution to (\ref{exact2}) 
can be found by approximating 
the exponential term in (\ref{exact2}) 
near $X=X^*$. 
The rate of change of the exponential's exponent near $X=X^*$ can be
estimated by considering the difference in its value between $X^*$ and
$(X^*-1)$, which is simply $\log[
(1+X^*/X_A)(1+X^*/X_B)/(1+X^*/N_f)]$. 
Provided this rate of change is small, then a Poisson distribution
will provide a good approximation near the maximum of (\ref{exact2}).   
If $X^*/X_A\ll 1$ and  $X^*/X_B\ll 1$ (implying $X^*/N_f\ll 1$), then 
$\log[ (1+X^*/X_A)(1+X^*/X_B)/(1+X^*/N_f)]$ will be small, and the
exponent will be approximately constant near $X^*$.   
Therefore if $X^*/X_A\ll 1$ and  $X^*/X_B\ll 1$, the Poisson
distribution provides a good approximation 
near the maximum of (\ref{exact2}).  
If a precise and accurate approximation for the moments of (\ref{exact2}) 
only requires a sufficiently 
precise approximation to (\ref{exact2}) near $X=X^*$ (we do not claim
to show this  
here), then the Poisson distribution 
will provide a good approximation for the moments of (\ref{exact2}).  
These remarks are consistent with the observations in the main text
that:  
(19) 
is always greater than 
(21), 
but provided
that $X^*/X_A\ll 1$ and  $X^*/X_B\ll 1$ (implying $X^*/N_f\ll 1$), the
exact 
(19) 
and approximated moments 
(21), 
are approximately the same (for a Poisson distribution, $\langle X
\rangle = X^*$ and $\sigma^2 = X^*$, {\sl e.g.} see 
Stirzaker (1994)). 
The above calculation easily generalises  to the case studied by
Garc\'ia-Pelayo (2006) 
with n-persons searching, consequently similar
remarks apply to that problem also. 



\section{Relation between (8) and (10)
}\label{newApp} 

Here the relationship between 
(8) 
and
(10)  
is discussed. 
Firstly write 
(8) 
in terms of $X=N-N_f$, $X_A = N_A
-N_{AB}$, $X_B = N_B - N_{AB}$, and $N_f=N_{AB}+X_A +X_B$, to give, 
\begin{equation}\label{ap1}
P(X|X_A,X_B,N_f) = 
\frac{(X+X_A)!(X+X_B)!}{X!(X+N_{f})!} 
\frac{P(X+N_f)}{(X+N_f+1)^2} C 
\end{equation}
Throughout this section we will only 
consider the case where $P(X+N_f)=P(N)$ is
constant. 
The moments of (\ref{ap1})  
are then, 
\begin{equation}\label{ap2}
\langle X^p \rangle = 
\frac{
\sum_{X=0}^{\infty} \frac{X^p}{(X+N_f+1)} 
\frac{(X+X_A)!(X+X_B)!}{X!(X+N_f+1)!} 
}
{
\sum_{X=0}^{\infty} \frac{1}{(X+N_f+1)} 
\frac{(X+X_A)!(X+X_B)!}{X!(X+N_f+1)!} 
}
\end{equation}
where one of the factors of $1/(X+N_f+1)$ has been incorporated into 
$1/(X+N_f+1)!$. 
Note that, 
\begin{equation}\label{s1}
\frac{1}{X+N_f+1} > \frac{1}{X+N_f+2} 
\end{equation}
and that, 
\begin{equation}\label{s2}
\begin{array}{ll}
\frac{1}{X+N_f+1} &= \frac{1}{X+N_f+2} 
\left( \frac{1}{1-\frac{1}{X+N_f+2}} \right)
\\
&= \frac{1}{X+N_f+2} \sum_{k=0}^{\infty} 
\left( \frac{1}{X+N_f+2} \right)^k 
\\
&< \frac{1}{X+N_f+2} \sum_{k=0}^{\infty} 
\left( \frac{1}{N_f+2} \right)^k = \frac{1}{X+N_f+2} 
\left(\frac{N_f+2}{N_f+1} \right) 
\end{array}
\end{equation}
Using these bounds (\ref{s1}) and (\ref{s2}) in the numerators and
denominators of (\ref{ap2}) as appropriate (with (\ref{s2}) used for
the sum in the numerator and (\ref{s1}) for the sum in the denominator
to give the upper bound, and vice versa for the lower bound), we find,  
\begin{equation}\label{c3}
\begin{array}{c}
\frac{1}{\left(\frac{N_f+2}{N_f+1}\right)} 
\frac{
\sum_{X=0}^{\infty} X^p 
\frac{(X+X_A)!(X+X_B)!}{X!(X+N_f+2)!} 
}
{
\sum_{X=0}^{\infty}  
\frac{(X+X_A)!(X+X_B)!}{X!(X+N_f+2)!} 
}
< \langle X^p \rangle 
\\
< 
\left( \frac{N_f+2}{N_f+1} \right) 
\frac{
\sum_{X=0}^{\infty} {X^p} 
\frac{(X+X_A)!(X+X_B)!}{X!(X+N_f+2)!} 
}
{
\sum_{X=0}^{\infty}  
\frac{(X+X_A)!(X+X_B)!}{X!(X+N_f+2)!} 
}
\end{array}
\end{equation}
where the factors of $1/(X+N_f+2)$ have been incorporated into the
factors of $1/(X+N_f+2)!$.  
Using $\langle X^p \rangle_0[N_f+2]$ to refer to moments of
(\ref{e30}), but with $N_f$ replaced by $N_f+2$, or equivalently
noting that $N_f = X_A + X_B + N_{AB}$, by replacing $N_{AB}$ by 
$N_{AB}+2$, keeping $X_A$ and $X_B$ fixed everywhere else.  
With this notation in (\ref{c3}), and using 
$-1/(N_f+1) < -1/(N_f+2)$ to make the left hand side of the inequality
symmetric with the right, we find,  
\begin{equation}
\left( 1 - \frac{1}{N_f+1} \right) 
\langle X^p \rangle_0 [N_f +2] < \langle X^p \rangle < 
\left( 1 + \frac{1}{N_f+1} \right) 
\langle X^p \rangle_0 [N_f +2]
\end{equation}
Or equivalently, 
\begin{equation}\label{a3-1}
\langle X^p \rangle = \langle X^p \rangle_0 [N_f +2] 
\left( 1 \pm \frac{1}{N_f+1} \right) 
\end{equation}
where the factor of $\pm 1/(N_f+1)$ gives a maximum error bound. 
Improved bounds can be found on a case by case basis, by considering 
$\langle X^p \rangle - \langle X^p \rangle_0[N_f+2]$, simplifying as
far as possible, then using (\ref{s1}) and (\ref{s2}) to express the
sums in a form that can be evaluated using (\ref{appS}).  
Returning to (\ref{a3-1}), if $N_f \gg 1$ then an excellent approximation to
$\langle X^p \rangle$ that is correct to within $\pm 100/(N_f+1)$ percent, is 
given by $\langle X^p \rangle_0[N_f+2]$.   
This approximation for the moments of (\ref{ap1}) 
is equal to the
exact moments of (\ref{e30}) 
with $N_{AB}$ replaced by
$N_{AB}+2$, keeping $X_A$ and $X_B$ fixed.  
Consequently using (\ref{avX})  
we have, 
\begin{equation}\label{avapp}
\langle X \rangle = 
\frac{(X_A+1)(X_B+1)}{N_{AB}} 
\mbox{ with } N_{AB} >0
\end{equation}
with a maximum error of $\pm \langle X \rangle/(N_f+1)$, which 
with the substitutions $X_A=N_A-N_{AB}$ and $X_B=N_B-N_{AB}$, is 
(23) 
of the main text. 
Similarly using (\ref{avX}) and (\ref{sigX}) 
an approximation for $\sigma^2$ is,  
\begin{equation}\label{relsig2}
\sigma^2 =   
\frac{(X_A+1)(X_B+1)(N_{AB}+X_A+1)(N_{AB}+X_B+1)} 
{N_{AB}^2\left(N_{AB}-1\right)} 
\mbox{ with } N_{AB} > 1
\end{equation}
which with the substitutions $X_A=N_A-N_{AB}$ and $X_B=N_B-N_{AB}$, 
is (24) 
of the main text. 
Unfortunately whereas (\ref{avapp}) 
has a maximum error of order
$\langle X \rangle/N_f$, which is much less than $\langle X \rangle$
if $N_f \gg 1$, $\sigma^2=\langle X^2 \rangle - \langle X \rangle^2$
can be arbitrarily small, but the maximum possible error remains of
order $\langle X^2 \rangle/N_f$.  
Therefore unless $\langle X^2 \rangle/N_f \ll 1$, (\ref{relsig2}) 
will not be  guaranteed to give a good approximation for $\sigma^2$.  
As is noted in the main text, an alternative approach is to use the
prior $P(N)=\kappa(N+1)/(N+2)$, for which 
(\ref{avapp}) and (\ref{relsig2}) 
are the exactly 
calculated moments. 
For that case this calculation gives the maximum
difference between the moments with this, and with a prior that is
independent of $N$. 

\end{document}